\documentclass[twoside, twocolumn]{article}

\usepackage{graphicx}
\usepackage{mathptmx}      
\usepackage{fancyhdr}
\usepackage[left = 2cm, right = 1.5cm]{geometry}
\usepackage{footmisc}
\usepackage{url}
\usepackage[numbers]{natbib}
\usepackage{amsmath}
\usepackage{amsfonts}
\usepackage{titlesec}

\titleformat*{\section}{\normalfont\bfseries}
\titleformat*{\subsection}{\normalfont}

\newcommand{\R}{\mathbb{R}}
\newcommand{\N}{\mathbb{N}}
\newcommand*\diff{\mathop{}\!\mathrm{d}}
\newcommand{\ueps}{u_\varepsilon}
\newcommand*{\defeq}{\mathrel{\vcenter{\baselineskip0.5ex \lineskiplimit0pt
                     \hbox{\scriptsize.}\hbox{\scriptsize.}}}%
                     =}
\newcommand*{\defeqrev}{=\mathrel{\vcenter{\baselineskip0.5ex \lineskiplimit0pt
                     \hbox{\scriptsize.}\hbox{\scriptsize.}}}%
                     }

\begin{document}

\title{Why more physics can help achieving better mathematics}

\author{Andr\'e Eikmeier\footnote{Institut f\"ur Mathematik, Technische Universit\"at Berlin, Stra\ss e des 17. Juni 136, 10623 Berlin, Germany, E-mail: \{eikmeier, emmrich\}@math.tu-berlin.de}      \and
        Etienne Emmrich\footnotemark[1]    \and
        Eckehard Sch\"oll\footnote{Institut f\"ur Theoretische Physik, Technische Universit\"at Berlin, Hardenbergstra\ss e 36, 10623 Berlin, Germany, E-mail: schoell@physik.tu-berlin.de}
}

\date{}

\maketitle

\begin{abstract}
In this paper, we discuss the question whether a physical ``simplification" of a model makes it always easier to study, at least from a mathematical and numerical point of view. To this end, we give different examples showing that these simplifications often lead to worse mathematical properties of the solution to the model. This may affect the existence and uniqueness of solutions as well as their numerical approximability and other qualitative properties. In the first part, we consider examples where the addition of a higher-order term or stochastic noise leads to better mathematical results, whereas in the second part, we focus on examples showing that also nonlocal models can often be seen as physically more exact models as they have a close connection to higher-order models.
\end{abstract}

\section{Introduction}
\label{intro}
Differential equations are nowadays a standard tool to describe many physical phenomena and processes. The goal of this paper is to illustrate that, in many cases, a physical ``simplification'' of a model (meaning here the omission of higher-order derivatives or stochastic noise) does not imply a mathematical simplification as it reduces the mathematical ``quality'' of solutions in terms of regularity, uniqueness, numerical approximability and other properties such as, for example, long-term behaviour. This raises the question if it might be more adequate to study the physically more exact model with respect to mathematics and numerics. Besides, we also wish to draw attention to the less common nonlocal models: In many cases, the equation with higher-order terms can be transformed into a nonlocal equation, or it is at least an approximation, in some sense, of a nonlocal equation. In this way, nonlocal models can be seen as models comprising different local models as well as models of different characteristic scale.\footnote{Sometimes, nonlocal models can be seen as an upscaling of a local, microscale model, for example the theory of peridynamics \cite{SPGL09}. The purpose of peridynamic models is to be computationally faster than molecular dynamics but to preserve characteristic properties of molecular dynamics that are not recovered by classical continuum mechanics.} This connection might also help improving the numerics of nonlocal equations, for which up to now less is known in comparison to the numerics of local models.

In the first part of this paper, we give a short introduction to some mathematical solution concepts for differential equation problems to make the following sections of this work easier to understand.

In the second part, we present some examples where the simplified model has a mathematically ``worse'' solution than the so-called regularised problem including higher-order terms or stochastic perturbations, which both have a physical meaning in most cases. We give a concrete example where the higher-order term already appears in the physical derivation of the model but is then omitted because it seems to be of negligibly small order of magnitude.

In the third and last part, we give some examples showing interesting connections between nonlocal models and local higher-order models, for example the transformation of a higher-order equation into a nonlocal one, or the approximation of a nonlocal equation via a higher-order equation.

\section{Different solution concepts} \label{solution concepts}

We briefly introduce some of the most important mathematical solution concepts for differential equation problems. Throughout this paper, let $\Omega$ be an open subset of $\R^d$ with sufficiently smooth boundary. We start with the probably most popular concept, the concept of classical solutions: A classical solution is sufficiently often continuously differentiable and the differential equation is fulfilled pointwise everywhere. A similar concept is the concept of strong solutions:\footnote{The term ``strong" solution is not uniquely determined in the mathematical literature. Many authors use it as a synonym for classical solutions, but we will stick here with the definition given in \textsc{Schweizer}~\cite{S13}, since in this way, it is clearly separated from the other definitions.} A strong solution is sufficiently often differentiable in the weak sense and the equation is fulfilled pointwise almost everywhere. We call a function \mbox{$u\colon \Omega\to \R$} weakly differentiable in the direction $x_i$, $i\in\{1,...,d\}$, if it, together with its weak derivative $v\colon \Omega\to \R$, fulfils some kind of integration-by-parts condition, i.e., the equation
\begin{equation} \label{weak derivative}
  \int_\Omega u(x) \;\! \partial_{x_i}\varphi(x) \diff x = - \int_\Omega v(x)\;\! \varphi(x) \diff x
\end{equation}
holds for all test functions $\varphi\colon \Omega\to \R$ that are infinitely many times continuously differentiable and vanish outside of a compact subset of the domain $\Omega$. This concept of weak derivatives can be seen as a generalisation of the concept of classical derivatives since every integrable classical derivative fulfils equation \eqref{weak derivative}. The boundary term appearing if integration by parts is applied vanishes here because the function $\varphi$ vanishes at the boundary of $\Omega$.

The weak derivative also motivates another quite popular concept of solutions, the so-called weak solutions. One of the first using it was \textsc{Jean Leray} \cite{L34} in 1934 to prove the existence of solutions to the \textsc{Navier--Stokes} equations in two and three spatial dimensions. Since then it has become one of the most often used concepts of solutions. To obtain the weak formulation of the differential equation problem, we multiply the original differential equation with test functions $\varphi$ with the same properties as in the definition of the weak derivative. Then we integrate on both sides of the equation over $\Omega$ and formally perform an integration by parts in the terms with higher derivatives such that half of the derivatives is shifted onto the test function $\varphi$. The solution to the resulting integral equation in an appropriate solution space we call weak solution. Note that the weak solution in general has a lower regularity than both classical and strong solutions, e.g., a weak solution to the common heat equation
\begin{equation}
  \partial_t u + \alpha\;\! \Delta u = f
\end{equation}
only has to be once differentiable in space in the weak sense although two derivatives appear in the original equation. It can be shown that a classical solution is also always a weak solution. Therefore, the concept of weak solutions is a generalisation of the concept of classical solutions. It allows for example discontinuities in some derivatives of the solution and is able to handle more general data. Moreover, there are many (especially nonlinear) problems that do not admit classical solutions but weak solutions. A more detailed introduction into the concept of weak solutions can, e.g., be found in \textsc{Chipot} \cite{C09} or \textsc{Roub\'i\v{c}ek} \cite{R13}.

The last concept of solutions we want to discuss in detail is the concept of the so-called \mbox{(\textsc{Young}-) measure}-valued solutions. In the case of a nonlinear equation, even the existence of weak solutions sometimes cannot be proven, especially if certain oscillation effects come into play. These oscillation effects are now covered by a measure-valued mapping $\nu$ which replaces the function $u$ (or derivatives of $u$, respectively) in the weak formulation in the term with the nonlinearity. Thus, a measure-valued solution always consists of the function $u$ and the measure-valued mapping $\nu$. The concept of measure-valued solutions is again a generalisation of the concept of weak solutions, since every weak solution is also a measure-valued solution taking as $\nu$ the point measure of the solution (or derivatives of it, respectively). However, the concept of measure-valued solutions is quite controversial, as it is a rather weak concept. For example, uniqueness of measure-valued solutions cannot be expected in general. Even for problems that admit a unique weak solution, existence of infinitely many measure-valued solutions can be shown \cite{RH95}. Apart from that, the numerical approximability is worse than for weak solutions, which, because of their nice structure, can be easily approximated via the finite element method, for example. A more detailed introduction into the concept of measure-valued solutions can, e.g., be found in \textsc{M\'alek} et al. \cite{MNRR96}.

Between these concepts there are many others, for example very weak solutions where, in contrast to weak solutions, all derivatives are shifted onto the test function, or entropy solutions which have to fulfil an additional entropy inequality. An overview of some of these concepts using the example of the \textsc{Navier--Stokes} equations can be found in \textsc{Amann}~\cite{A00}. Another overview of results concerning classical, weak and measure-valued solutions in the theory of elastodynamics is given in \textsc{Emmrich} and \textsc{Puhst}~\cite{EP15}. 

\section{Higher-order terms in physics and mathematics} \label{higher-order}

In this section, we present two examples where physically motivated higher-order terms lead to mathematically ``better'' solutions to the regularised problem than to the original problem. Additionally, we consider the regularisation using stochastic perturbations instead of higher-order terms, which requires less assumptions on the given data, compared to the problem without a stochastic perturbation, to obtain existence and uniqueness of solutions.

\subsection{Example 1: Backward-forward heat equation}\label{bfhe}
  
  The first example is the so-called backward-forward heat equation
  \begin{equation} \label{bf_heat_eq}
    \partial_t u - \nabla\cdot\phi(\nabla u)=0,
  \end{equation}
  where $u$ represents the temperature and $\phi$ is a nonmonotone function representing the heat flux density. This equation occurs, for example, in cases where \textsc{Fourier}'s law of heat conduction cannot be applied to simplify the general heat equation. Apart from thermodynamics, it is also very important for the so-called anisotropic diffusion considered in image processing.
  
  There are several articles proving existence of measure-valued solutions to the back\-ward-forward heat equation, for example, \textsc{Slemrod} \cite{S91} and \textsc{Thanh} et al. \cite{TST14}. In both articles, the method of regularisation is used to prove existence. In the first work \cite{S91}, the term $\varepsilon\;\! \Delta^2 u$ is added to the left-hand side, where $\varepsilon>0$ is small. For the regularised equation 
  \begin{equation}
    \partial_t u_\varepsilon - \nabla\cdot\phi(\nabla u_\varepsilon) + \varepsilon\;\! \Delta^2 u_\varepsilon	= 0,
\label{heat_regularized}
  \end{equation}
  existence of weak solutions can be proven. The limit \mbox{$\varepsilon\to0$} then yields the existence of a measure-valued solution to Eq.~\eqref{bf_heat_eq}. As mentioned in Section \ref{solution concepts} on solution concepts, measure-valued solutions are weaker than weak solutions. Thus, we get ``better'' solutions to the regularised equation with a higher-order term than to the original equation without this higher-order term.
  
  In the second work \cite{TST14}, the term $-\varepsilon\;\!\Delta \partial_t u$ is added to the left-hand side of the backward-forward heat equation. Again, existence of weak solutions can be shown for the regularised equation 
  \begin{equation}
    \partial_t \ueps-\nabla\cdot\phi(\nabla u_\varepsilon) - \varepsilon\;\!\Delta \partial_t \ueps =0,
  \end{equation}
  but in the limit $\varepsilon\to 0$, we again only end up with a measure-valued solution to the original equation.
  
  So far, the existence of weak solutions to the backward-forward heat equation is only known in a special case, where the spatial dimension is equal to one and where the heat flux density $\phi$ is piecewise linear (cf. \cite{H83,HN83}). Besides, uniqueness of solutions is also a problem. Uniqueness of measure-valued solutions cannot be expected due to the solution concept, and in the case mentioned above, where existence of weak solutions is known, it can additionally be proven that there exist infinitely many weak solutions. So far, uniqueness is only known for a special kind of classical solutions and again only in the case of one spatial dimension and for special cases of the heat flux density $\phi$ (cf. \cite{L85,L88}). Whether solutions of this kind do even exist, could, to the best knowledge of the authors, not yet be proven.
  
  These observations raise one central question of this article: Is passing to the limit $\varepsilon\to 0$ and thus reducing the mathematical quality of solutions neces\-sary for the equations to be an adequate model of the underlying physics? In fact, both regularisations shown above are motivated by physics. \textsc{Slemrod} \cite{S91} refers to the higher-order theory of heat conduction due to \textsc{Maxwell} \cite{M79}, which is also mentioned in \textsc{Truesdell} and \textsc{Noll} \cite{TN65}. It is based upon a moment expansion of the \textsc{Boltzmann} equation, which leads to an infinite hierarchy of coupled moment equations for the mass density, the momentum density, the energy density, the energy flux etc., and has to be truncated by appropriate closure approximations. The common approximation of the heat flux by \textsc{Fourier}'s law in the mean energy balance equation gives the standard heat conduction equation. However, already in \cite{M79} it was shown for rarified gases that higher-order terms can arise due to density gradient contributions to the energy balance, and in \cite{S91} it was pointed out that a double Laplacian similar to Eq.~\ref{heat_regularized} also occurs in the \textsc{Cahn}--\textsc{Hilliard} equation which describes the process of phase separation in a two-component binary fluid. Moment expansions of the \textsc{Boltzmann} equation, resulting in hydrodynamic balance equations for charge carrier densities and electron temperature, have also widely been used in nonlinear electron transport in semiconductors \cite{SCH98,QUA91,QUA93a}, giving higher-order terms at various levels of approximation.
  
  \textsc{Thanh} et al. \cite{TST14} describe the regularisation term as some kind of viscous effects referring, amongst others, to \textsc{Binder}~et~al.~\cite{BFJ86} and \textsc{Novick-Cohen} and \textsc{Pego} \cite{NCP91}. 
  
  Indeed, there is a whole mathematical theory called \textit{vanishing viscosity method} based upon the regularisation with a diffusive term multiplied with the regularisation parameter~$\varepsilon$. At the end, letting the parameter~$\varepsilon$ tend to zero yields a solution to the original problem. The method goes back to \textsc{von Neumann} and \textsc{Richtmyer} \cite{NR50} who introduced it as a method to calculate hydrodynamic shocks numerically.
  
\begin{figure}
  \includegraphics[width=\columnwidth]{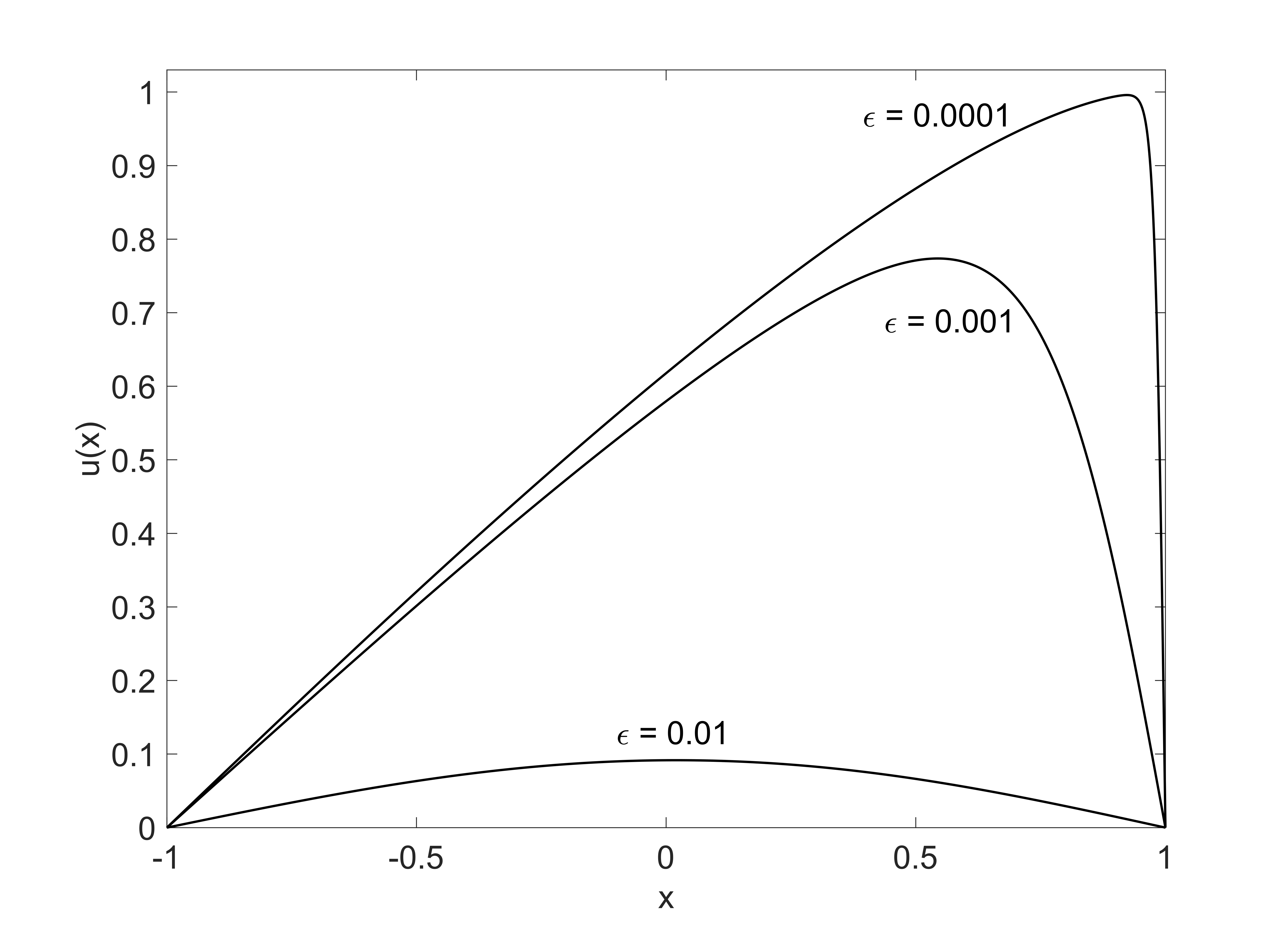}
\caption{Solutions to the \textsc{Burger}'s equation at a finite time equipped with the initial condition $u_0(x)=-x^2+1$ and homogeneous \textsc{Dirichlet} boundary conditions for different values of $\varepsilon$}
\label{fig:burgers}       
\end{figure}

  One of the simplest examples for which the \textit{vanishing viscosity method} can successfully be applied is the famous \textsc{Burgers}' equation
  \begin{equation} \label{Burgers}
    \partial_t u + u \;\!\partial_x u =0,
  \end{equation}
  which can be seen as a simple model for a one-dimensional flow. To apply the \textit{vanishing viscosity method}, the diffusion term $\varepsilon\;\! \partial_{xx}u$ is added,
  \begin{equation}
    \partial_t u + u \;\!\partial_x u -\varepsilon\;\! \partial_{xx}u=0,
  \end{equation}  
  where $\varepsilon$ is the regularisation parameter mentioned above. Intuitively, the discontinuities that may appear in the case of the inviscid \textsc{Burgers}' equation \eqref{Burgers} are first smoothed out by the additional diffusion term and then, this smoothing effect is reduced more and more (cf. Figure \ref{fig:burgers}), such that at the end, a solution to the original equation \eqref{Burgers} is obtained. A detailed discussion of the application of the \textit{vanishing viscosity method} in the case of the \textsc{Burgers}' equation can, e.g., be found in the lecture notes of \textsc{Kruzhkov} on first-order quasilinear partial differential equations \cite{CG09}.
  
   Another example where the \textit{vanishing viscosity method} can successfully be applied is a system of \textsc{Boussinesq} equations, shown in the monograph of \textsc{Guo} et al. \cite{GBLX17}. It reads
  \begin{equation}
    \begin{aligned}
      \partial_t \rho + \alpha\;\! \partial_x u + \beta\;\! \partial_x (u\rho)& =0, \\
      \partial_t u + \gamma\;\!\partial_x \rho + \delta\;\! u\;\! \partial_x u - \nu\;\! \partial_{xxt} u &=0,
    \end{aligned}
  \end{equation}
  and describes the propagation of the long surface wave in a pipe with constant depth. Here, $\rho$ is the density, $u$ is the velocity, $\omega=1+\delta\rho$ denotes the altitude from the bottom to the free surface of the flow and $\alpha,\beta,\gamma,\delta, \nu$ are constants. In this case, the regularisation term $\varepsilon\;\! \partial_{xx} \rho$ is added to the first equation and at the end, $\varepsilon$ tends to zero. Here, the regularised problem admits a unique classical smooth solution, whereas the original problem only admits weak solutions.
  
  A broad overview of other examples for the \textit{vanishing viscosity method} can also be found in \textsc{Guo} et al. \cite{GBLX17}.
  
\subsection{Example 2: Kinetic models for dilute polymers} \label{dilute polymers}
  
  In other examples, like the following one, the regularisation term is already existent in the derivation of the physical model but is then omitted because it seems to be of negligibly small order of magnitude. \textsc{Barrett} and \textsc{S\"uli} \cite{BS07} consider a kinetic bead-spring model for dilute polymers, where the extra-stress tensor is defined through the associated probability density function $\psi$. This function satisfies the \textsc{Fokker--Planck}-type parabolic equation
  \begin{multline} \label{fokpla}
    \partial_t \psi + \left( u\cdot \nabla_x\right)\psi+ \nabla_q\cdot \left(\left( \nabla_x\;\! \mathcal{J}^x_{l_0,q}\;\! u \right)q\;\! \psi\right) \\
    = \varepsilon\;\! \Delta_x \psi+ \frac{1}{2\lambda}\nabla_q \cdot \left(\nabla_q\;\! \psi + U'\;\!q\;\! \psi\right)\!.
  \end{multline}
  Here, $u$ is the velocity of the fluid considered, $q$ is the elongation vector of the dumbbell representing a polymer chain, $\mathcal{J}^x_{l_0,q}$ is the directional \textsc{Friedrichs} mollifier with respect to $x$ over an interval of length $l_0\vert q\vert$ in the direction $q$, and $U$ is the potential of the elastic force of the spring connecting two beads. The constant $\varepsilon$ corresponds to the quantity $\frac{\text{De}}{\text{Pe}}$, where $\text{De}$ denotes the \textsc{Deborah} number and $\text{Pe}$ the \textsc{P\'eclet} number, and the constant $\lambda$ corresponds to the relaxation time constant of the dumbbells.
  
  One interesting feature of this model is the presence of the diffusion term $\varepsilon\;\!\Delta_x\;\!\psi$ on the right-hand side of Eq.~\eqref{fokpla}. As \textsc{Barrett} and \textsc{S\"uli} \cite{BS07} already mention, this term is usually omitted in standard derivations of bead-spring models, because it is several orders of magnitude smaller than the other terms in Eq. \eqref{fokpla}. Actually, \textsc{Bhave}~et~al.~\cite{BAB91} estimate the quantity $\frac{De}{Pe}$ to be in the range of $10^{-9}$ to $10^{-7}$, whereas the expected important length scales of stress diffusion start at $10^{-5}$ to $10^{-3}$.
  
    Mathematically however, omitting the diffusion term is quite detrimental as it leads to a hyperbolically degenerate parabolic equation which is much harder to handle than Eq.~\eqref{fokpla}.\footnote{If a parameter tends to zero and thus changes the type of the equation, this phenomenon is called singular limit. Nice overviews of such singular limits in hydrodynamics and thermodynamics of viscous fluids can be found in \textsc{Masmoudi} \cite{M07} and \textsc{Feireisl} and \textsc{Novotn\'y}~\cite{FN09}, respectively.} In fact, the existence result in the case $\varepsilon=0$ is again proven via showing the existence of solutions for $\varepsilon>0$ and then passing to the limit as~$\varepsilon\to 0$. This leads to less regularity for the probability density function $\psi$. So again we state that the original model with the higher-order term delivers mathematically ``better'' solutions than the model without this higher-order term. In this case, the model with the higher-order term even seems to be more adequate physically.
  
\subsection{Example 3: Regularisation by noise} \label{noise}

A quite recent approach to regularise an equation is to add a certain stochastic noise in order to obtain the existence of a unique solution where, without noise, only existence or uniqueness or none of these two is known so far. 

There is some work by, e.g., \textsc{Gy\"ongy} and \textsc{Pardoux} \cite{GP93a,GP93b} using additive noise to prove existence of a unique solution under assumptions which, in the deterministic case, are so far not known to suffice for obtaining existence or uniqueness. \textsc{Gy\"ongy} and \textsc{Pardoux} \cite{GP93a,GP93b} consider the equation
\begin{equation}
  \partial_t u(x,t) - \partial_{xx} u(x,t) =f(x,t,u(x,t)) + \partial_{tx} W(x,t)
\end{equation}
equipped with either homogeneous \textsc{Dirichlet} or homogeneous \textsc{Neumann} boundary conditions, where $f$ is a nonlinear function and $\partial_{tx} W$ denotes space-time white noise. Under the assumption that $f$ satisfies some measurability and boundedness condition, \textsc{Gy\"ongy} and \textsc{Pardoux} \cite{GP93a,GP93b} are able to prove existence and uniqueness of a solution in a generalized sense defined in the work of \textsc{Walsh} \cite{W86}, which may be compared to a very weak solution (cf. Section \ref{solution concepts}) but in a stochastic sense. In the deterministic case, the assumptions on $f$ are, to the best knowledge of the authors, not enough to prove existence or uniqueness of solutions.

More recently, there has been research on linear multiplicative noise by, e.g., \textsc{Flandoli}~et~al.~\cite{FGP10}, considering the linear transport equation
\begin{equation} \label{linear transport}
  \partial_t u + b\cdot \nabla u = 0
\end{equation}
driven by the vector field $b$. Assuming that $b$ is sufficiently regular, uniqueness of solutions to the initial-boundary value problem governed by this equation can be proven (cf., e.g., \textsc{DiPerna} and \textsc{Lions} \cite{DL89} or \textsc{Ambrosio} \cite{A04}), but if this is not the case then examples of non-uniqueness are known, as is shown in the work of \textsc{Flandoli}~et~al.~\cite{FGP10}. However, if a certain amount of linear multiplicative noise is added to Eq. \eqref{linear transport}, existence and uniqueness of solutions can be proven under weaker assumptions on $b$, see again \cite{FGP10}. To be precise, the stochastic equation
\begin{equation}
  \diff_t u + (b\cdot \nabla u)\diff t + \sum_{i=1}^d e_i\cdot \nabla u \circ \diff W^i_t = 0
\end{equation}
is considered, where $e_i$, $i=1,...,d$, are the unit vectors in $\R^d$, $W_t\defeq (W_t^1,..., W_t^d)$ is a standard \textsc{Brownian} motion in $\R^d$, and the notation $\circ$ is used for the stochastic integration in the sense of \textsc{Stratonovich}.

Since real-world systems often include noise, the consideration of stochastic differential equations is physically also very important, and there are many works considering the influence of noise on various physical systems. We just want to mention some examples here. Additive noise has been shown to have an important effect, e.g., upon chimera states (coexisting coherent and incoherent space-time patterns in networks), which can be either destructive, see, e.g., \textsc{Loos} et al. \cite{LZCS15,ZAK16}, or constructive, see, e.g., \textsc{Zakharova} et al. \cite{SEM16,ZAK17}. Multiplicative noise has been considered, for example, in the work on nonequilibrium phase transitions by \textsc{Van den Broeck} et~al.~\cite{BPTK97}.  
  
\section{Connection of higher-order and nonlocal equations} \label{nonlocal}
  
  In this part, we show some interesting connections between nonlocal models and local higher-order models. As it turns out, many equations containing higher-order terms can be rewritten as some nonlocal equation or can be seen as approximations of nonlocal equations. Thus, nonlocal models can also be seen as physically more exact regularisations of local models.
  
\subsection{{\normalfont\textsc{Green}}'s function as integral kernel} \label{Green's function}
  
  Let us start with a simple example, which was, amongst others, mentioned in \textsc{Duruk}~et~al.~\cite{DEE09,DEE10,DEE11}. We consider the equation
  \begin{equation} \label{nonlocal wave}
    \partial_{tt} u(x,t)- \int_{\Omega} G(x-\xi)\;\! \partial_{xx} u(\xi,t)\diff \xi = 0
  \end{equation}
  arising, e.g., in the nonlocal theory of elasticity which considers the stress at a point $x$ not only as a function of the strain at the point $x$ but the strain field at every point in the body. Thus, it also takes long-range interactions between the particles into account. In this case, $u$ represents the displacement field and $G$ is the kernel function representing the stress-strain relation. Here, $G$ is the \textsc{Green}'s function for the differential operator $1-\varepsilon\;\!\partial_{xx}$ considered on the full space (therefore we can neglect boundary conditions), i.e.,
  \begin{equation}
    v(x)\defeq \int_\Omega G(x-\xi)\;\!f(\xi)\diff\xi
  \end{equation}
  is the unique solution to the problem
  \begin{equation}
	  v(x) - \varepsilon\;\!\partial_{xx} v(x) = f(x).
  \end{equation}
  Based on Eq. \eqref{nonlocal wave}, we can now conclude that $\partial_{tt} u$ at least formally fulfils
  \begin{equation}
	      \partial_{tt} u(x,t) - \varepsilon\;\!\partial_{xx} \partial_{tt} u(x,t) = \partial_{xx} u(x,t).
  \end{equation}
  Therefore, we can see the nonlocal equation \eqref{nonlocal wave} as a regularised form of the usual wave equation
  \begin{equation}
    \partial_{tt} u(x,t) - \partial_{xx} u(x,t)=0.
  \end{equation}
  
  This method can obviously be applied using any arbitrary \textsc{Green}'s function $G$ in an arbitrarily high spatial dimension, generating many different regularisations of arbitrarily high order. Considering, for example, the general differential operator
  \begin{equation}
    A\defeq  \mathop{Id}+\sum_{k=1}^n (-1)^k a_k\;\!\Delta^k,
  \end{equation}
  where $n\in\N$, $a_k>0$, $k=1,...,n$, the \textsc{Green}'s function for $A$ reads
  \begin{equation}
    G(x)= \mathcal{F}^{-1}\left(\left(1+\sum_{k=1}^n a_k \vert \omega \vert^{2k}\right)^{-1}\right),
  \end{equation}
  where $\mathcal{F}^{-1}$ denotes the inverse \textsc{Fourier} transform. If we start, for example, with a nonlocal equation, this operator $A$ can be used to approximate the differential operator of which the given kernel function is the \textsc{Green}'s function. This has been done, for example, for lattice models in nonlocal continuum mechanics in \textsc{Eringen} \cite[Ch. 6.9]{E02} and \textsc{Lazar} et~al.~\cite{LMA06}. Thus, we can consider nonlocal models as a generalisation of local ones as different kernels lead to different local, potentially higher-order models.
  
\subsection{Nonlocal coupling in reaction-diffusion systems} \label{coupling}
  
  A more involved example is a reaction-diffusion system with asymmetric nonlocal coupling, studied in \textsc{Siebert} et al. \cite{SABS14}, which is derived as a limiting case of the activator-inhibitor reaction-diffusion-advection model
  \begin{equation} \label{reaction-diffusion}
    \begin{aligned}
      \partial_t u(x,t) &= F(u(x,t)) - g\;\!w(x,t) + \partial_{xx} u(x,t), \\
      \tau \partial_t w(x,t) &= h\;\!u(x,t) - f\;\! w(x,t)+ \xi\;\! \partial_x w(x,t) +D\;\! \partial_{xx} w(x,t).
    \end{aligned}
  \end{equation}    
  Here, $u$ and $w$ are the activator and the inhibitor, respectively, linearly coupled by the terms $-g\;\!w$ and $h\;\!u$ with $g,h\in\R$, $F$ is a nonlinear function, $\tau$ is the inhibitor relaxation time, $f$ is a constant, $\xi$ represents the advection velocity and $D>0$ is the inhibitor diffusion coefficient. Passing to the limit \mbox{$\tau\to 0$}, i.e., the case of fast inhibitor dynamics, the \textsc{Green}'s function is used to reformulate the system by writing the solution $w$ of the second equation of \eqref{reaction-diffusion} as
  \begin{equation} \label{nonlocal_w}
    w(x,t) = h \int_{-\infty}^\infty G(y) \;\! u(x-y,t) \diff y
  \end{equation}
  with the \textsc{Green}'s function $G$ defined as
  \begin{equation}
    G(x) = \frac{f}{\sqrt{\xi^2 + 4Df}} e^{-(\sqrt{\xi^2 + 4Df}/2D)\;\!\vert x\vert} e^{\xi x/2D}.
  \end{equation}
  Inserting Eq. \eqref{nonlocal_w} into the first equation of \eqref{reaction-diffusion}, the reaction-diffusion equation with asymmetric nonlocal coupling
  \begin{equation}
    \partial_t u = F(u) + \partial_{xx}u - \sigma \int_{-\infty}^\infty G(y) \;\! u(x-y,t) \diff y
  \end{equation}
  is obtained, where $\sigma = g\;\!h$ denotes the nonlocal coupling strength.
  
  This method is more or less the inverse of the method shown in Section \ref{Green's function}. If we solve the first equation of \eqref{reaction-diffusion} for $w$ and insert it into the second equation, we end up with the same local higher-order equation which we would get using the method of Section \ref{Green's function}. Thus, higher-order and nonlocal terms may be considered as physically representing some kind of nonlocal coupling.
  
\subsection{Expansion of integral operators} \label{expansion integral}
  
  Another way to approximate a nonlocal equation by a local one including higher-order terms is the expansion of the corresponding integral operator. We will illustrate this by the example of a linear peridynamic model,\footnote{The peridynamic model goes back to \textsc{Silling} \cite{S00}, who introduced it as an alternative to the classical elasticity model, replacing the local stress by a nonlocal one independent of the spatial gradients of deformation.} as was done in \textsc{Emmrich} and \textsc{Weckner} \cite{EW07}. We consider the equation
  \begin{equation} \label{peridynamic nonlocal}
    \begin{aligned}
      &\rho(x)\;\! \partial_{tt} u(x,t) \\
      &= \int_{\Omega\cap B(x;\delta)} \lambda_{d,\delta}(\vert\hat{x}-x\vert) (\hat{x}-x) \otimes (\hat{x}-x) \\
      &\hspace{3cm}\cdot \left( u(\hat{x},t) - u(x,t) \right) \diff\hat{x} + b(x,t) \\
      &\defeqrev (L_{d,\delta} u)(x,t) + b(x,t),
    \end{aligned}
  \end{equation}
  where $\rho$ denotes the mass density, and $u$ is the displacement field of the body, $\Omega$ represents the volume that the body occupies, $B(x,\delta)$ is the ball of radius $\delta$ around $x$, $d$ is the spatial dimension, $\delta>0$ is the so-called peridynamic horizon of interaction, $\lambda_{d,\delta}$ is a real-valued function which depends upon $d$ and $\delta$ and determines the specific material model, and $b$ represents the external forces. The symbol $\otimes$ denotes the outer product. Additionally, the function $\lambda_{d,\delta}$ has the property
  \begin{equation}
    \lambda_{d, \delta}(r) = 0 ~~~ \text{for all } r \geq \delta.
  \end{equation}
  
  We now consider only the integral operator $L_{d,\delta}$. The idea is to perform a \textsc{Taylor} expansion of the argument $u$ up to the order of $m$, supposing sufficient regularity. In order to make the following easier to read, we omit the time dependence of $u$ and obtain
  \begin{equation}
    u(\hat{x}) - u(x) = \sum_{k=1}^m \frac{1}{k!} \left( (\hat{x}-x)\cdot \nabla\right)^k u(x) + r_m(u; \hat{x}, x) 
  \end{equation}
  with some remainder $r_m(u;\hat{x},x)$ of order $o (\vert \hat{x} - x \vert^m)$.
  Now we insert the \textsc{Taylor} expansion into the operator $L_{d,\delta}$ and get
  \begin{equation}
    (L_{d,\delta} u)(x) = \sum_{k=1}^m (L_{d,\delta}^{(k)} u)(x) + R_{m;d,\delta}(u;x)
  \end{equation}
  with
  \begin{multline} \label{L^k}
    (L_{d,\delta}^{(k)} u)(x) \defeq \frac{1}{k!} \int_{B(x;\delta)} \lambda_{d,\delta} (\vert\hat{x}-x\vert) (\hat{x}-x) \otimes (\hat{x}-x)\\
     \left( (\hat{x} -x) \cdot \nabla\right)^k u(x) \diff\hat{x}
  \end{multline}
  and again some remainder $R_{m;d,\delta}(u;x)$ of order $o(\delta^{m-2})$ in an appropriate norm (for details see \textsc{Emmrich} and \textsc{Weckner} \cite[Theorem 3.2]{EW07}). Since all derivatives appearing in the integrand of Eq.~\eqref{L^k} only depend on $x$ and not on $\hat{x}$, the operators $L_{d,\delta}^{(k)}$ are local differential operators of order $k$. Furthermore, for odd $k$, the integrand is an odd function in $\hat{x}-x$ so that the integral vanishes and thus $(L_{d,\delta}^{(k)} u)(x) =0$ for all odd $k$. In summary, we can approximate the nonlocal equation \eqref{peridynamic nonlocal} by the local equation
  \begin{equation}
    \rho(x)\;\! \partial_{tt} u(x,t) = \sum_{k=1}^{[m/2]} (L_{d,\delta}^{(2k)} u)(x) + b(x,t).
  \end{equation}
  For the second order $m=2$, this approximation yields the classical \textsc{Navier} equation of linear elasticity
  \begin{equation}
    \rho(x)\;\! \partial_{tt} u(x,t) = \mu \Delta u(x,t) +(\lambda+\mu) \nabla (\nabla \cdot u(x,t))+ b(x,t)
  \end{equation}
  with the \textsc{Lam\'e} coefficients $\mu=\lambda=3K/5$ where $K$ denotes the bulk modulus. Using the above approximation, it is shown that the nonlocal operator $L_{d,\delta}$ converges to the local operator $L$ of the \textsc{Navier} equation for vanishing nonlocality $\delta\to 0$, at least in the interior of the domain~$\Omega$. Near the boundary, this is to the best knowledge of the authors not yet known. The reason for that is the problem of finding appropriate boundary conditions for the \textsc{Navier} equation. For the peridynamic model, there are no boundary conditions since no spatial derivatives appear. More information about vanishing nonlocality in linear and also nonlinear peridynamic models can be found, e.g., in \textsc{Puhst} \cite{P16}.
  
  A quite similar approach to the one above, the so-called inner expansion, was used by \textsc{Arndt} and \textsc{Griebel}~\cite{AG05} to derive a new scheme upscaling the atomistic level to the continuum level for the model of a crystalline solid. There, a \textsc{Taylor} expansion of the deformation function around the points of the discrete lattice representing the atoms is inserted into the local energy potential to obtain, after some more steps, a continuum formulation which represents the upscaled model. One of this technique's advantages is the well-posedness of the resulting equation, which is not guaranteed by other techniques like the direct expansion technique proposed by \textsc{Kruskal} and \textsc{Zabusky} \cite{KZ64} and \textsc{Za\-busky} and \textsc{Kruskal}~\cite{ZK65}, which relies on a \textsc{Taylor} expansion of the right-hand side of the equation instead of the deformation function. On the other hand, it covers higher-order effects to an arbitrarily high order in contrast to the common scaling technique, which, simply spoken, lets the number of atoms tend to infinity, and was studied, amongst others, by \textsc{Blanc}~et~al.~\cite{BLL02}. This inclusion of higher-order effects prevents the loss of regularity of the solution that is often observed for the scaling technique because the dispersion which the discrete system inherits is not contained in the resulting continuum system.
  
  Again we see how important higher-order terms are and, although there is no nonlocal model involved in this example, there is still a connection to nonlocal models since they can also be seen as some kind of upscaling of a microscale model, as already mentioned in the introduction.
  
\section{Conclusion} \label{conclusion}

In the first part of this paper, we presented many different examples showing that simplifications of differential equations in the physical sense often lead to worse mathematical properties of the solutions to these equations. This may affect the existence of various types of solutions, the uniqueness of these solutions if they exist, their numerical approximability, and their qualitative properties such as long-term behaviour. The example of the kinetic model for dilute polymers (Section~\ref{dilute polymers}), studied in \textsc{Barrett} and \textsc{S\"uli} \cite{BS07}, shows that the higher-order term -- artificially added to the equation in order to regularise it -- is in some cases already existent in the original physical model. However, it is usually omitted because its order of magnitude is several orders smaller than the other terms appearing in the model. This sheds new light on the question if physical simplification is admissible, or if, on the contrary, even small terms should be kept in physical models since they improve the numerical solvability and the mathematical properties of the solution.

In the second part, we extended our results by demonstrating an intriguing connection between nonlocal and local higher-order models. The example of \textsc{Green}'s function as the integral kernel of the nonlocal equation (Section~\ref{Green's function}) shows that many local higher-order equations can be rewritten as a nonlocal equation, and in this sense nonlocal models can be seen as a generalisation of local models. The example of the expansion of integral operators (Section~\ref{expansion integral}) shows that even in more general situations, there is still a connection between nonlocal models and local higher-order models since the latter can be used as an approximation of the former. Overall, we wish to motivate the reader to consider also nonlocal models when studying a physical problem, especially in applications where nonlocal models have not received much attention previously. The connection mentioned might also improve the numerical understanding of nonlocal models.

As a physical interpretation of this connection, we suggest that nonlocal coupling terms can in certain situations be avoided by considering a more detailed description using additional dynamical variables. Considering a single differential equation with higher-order derivatives is equivalent to a system of several differential equations with lower-order derivatives, i.e., to extending the number of degrees of freedom of a dynamical model. These additional differential equations can be eliminated using \textsc{Green}'s functions, which results in nonlocal models with fewer dynamical variables. However, some nonlocal models like, for example, peridynamics would, at least formally, require an infinite number of degrees of freedom to be exactly represented. Additionally, they need much less regularity and are therefore more suitable to study discontinuous phenomena like, e.g., crack propagation.

\section*{Acknowledgements}
This work was supported by DFG in the framework of SFB 910, projects A1 and A8.

\bibliographystyle{mybibstyle}
{\footnotesize\bibliography{mybib}}   

\end{document}